\documentclass[11pt]{article}
\usepackage{latexsym}

\newcommand{\qed}{\ \hfill\mbox{$\Diamond$}\vspace{\baselineskip}}

\newtheorem{theorem}{Theorem}[section]
\newtheorem{conjecture}[theorem]{Conjecture}
\newtheorem{lemma}[theorem]{Lemma}
\newtheorem{proposition}[theorem]{Proposition}
\newtheorem{definition}{Definition}
\newtheorem{corollary}[theorem]{Corollary}

\newenvironment{proof}{\noindent {\bf Proof:}}{{\qed}}

\newcommand{\df}{:=}
\newcommand{\Nnn}{\mbox{\bf N}}
\newcommand{\HD}{\mbox{$D$}}
\newcommand{\iv}{\ensuremath{{\cal I}}}
\newcommand{\lv}{\ensuremath{\ell}}
\newcommand{\LV}{\mbox{$L$}}
\newcommand{\rank}{\rho}
\newcommand{\1}{\hat{1}}
\newcommand{\0}{\hat{0}}

\begin{document}

\title{Flag vectors of Eulerian partially ordered sets}

\author{Margaret M. Bayer\thanks{This research was supported by University of
Kansas General Research allocation \#3552.}\ \ and 
	G\'abor Hetyei\thanks{On leave from the
        Alfr\'{e}d R\'{e}nyi Institute of Mathematics,
        Hungarian Academy of Sciences.
	Partially supported by Hungarian National Foundation for
	Scientific Research grant no. F 023436.}\\
	Department of Mathematics\\
	University of Kansas\\
	Lawrence KS 66045-2142\\
}

\date{}

\maketitle

\begin{abstract}
The closed cone of flag vectors of Eulerian partially ordered sets is 
studied.
It is completely determined up through rank seven.
Half-Eulerian posets are defined.
Certain limit posets of Billera and Hetyei are half-Eulerian; they give
rise to extreme rays of the cone for Eulerian posets.
A new family of linear inequalities valid for flag vectors of Eulerian
posets is given.
\end{abstract}

\section{Introduction}

The study of Eulerian partially ordered sets (posets) originated with Stanley
(\cite{Stanley-groups}).
Examples of Eulerian posets are the posets of faces of regular CW spheres.
These include face lattices of convex polytopes, the Bruhat order on 
finite Coxeter groups, and the lattices of regions of oriented matroids.
(See \cite{Bjorner-CW} and \cite{BLVSWZ}.)

The flag $f$-vector (or simply flag vector) of a poset is a standard parameter 
counting chains in the partially ordered set by ranks.
In the last twenty years there has grown a body of work on numerical conditions
on flag vectors of posets and complexes, especially those arising in geometric
contexts.
Early contributions are from Stanley on balanced Cohen-Macaulay complexes
(\cite{Stanley-balanced}) and Bayer and Billera on the linear equations on
flag vectors of Eulerian posets (\cite{Bayer-Billera}).
A major recent contribution is the determination of the closed
cone of flag vectors of all graded posets by Billera and Hetyei
(\cite{Billera-Hetyei}).
Results on flag vectors and other invariants of Eulerian posets and special
classes of them are surveyed in \cite{Stanley-Eulerian}.

Our goal has been to describe the closed cone 
${\cal C}^{n+1}_{{\cal E}}$ of flag $f$-vectors 
of Eulerian partially ordered sets.
This problem was posed explicitly in~\cite{Billera-Liu}.
The ideal description would give explicitly both the facets (i.e., crucial
inequalities on flag vectors) and posets that generate the extreme rays. 
We have a complete solution only for rank at most seven.
For arbitrary ranks we give some of the facets and extreme rays.
The extreme rays of the general graded cone (\cite{Billera-Hetyei})
play an important role.
We introduce half-Eulerian partially ordered sets in order to incorporate
these limit posets in this work.

The remainder of this section provides definitions and other background,
and the definition of the flag $\LV$-vector, which simplifies the 
calculations.
Section~\ref{half-Eul} describes the extreme rays of the general graded
cone, defines half-Eulerian posets,
identifies which limit posets are half-Eulerian, and computes the corresponding
$cd$-indices.
Section~\ref{sec-inequalities} gives two general classes of inequalities on 
Eulerian flag vectors.
Section~\ref{sec-cone} shows that the half-Eulerian limit posets all give 
extremes of
the Eulerian cone, identifies some inequalities in all ranks as facet-inducing,
and describes completely the cone for rank at most 7.

\subsection{Background}

\label{ss_gpo}
A {\em graded poset} $P$ is a finite 
partially ordered set with a unique minimum element
$\0$, a unique maximum element $\1$, and a {\em rank function} $\rank:
P\longrightarrow \Nnn$ satisfying $\rank (\0)=0$, and
$\rank(y)-\rank(x)=1$ whenever $y\in P$ covers $x\in P$.
The {\em rank $\rank(P)$ of a graded poset $P$} is the rank of its maximum 
element. Given a graded poset $P$
of rank $n+1$ and a subset $S$ of $\{1,2,\ldots,n\}$ (which we abbreviate
as $[1,n]$), define the {\em
$S$--rank--selected subposet of $P$} to be the poset
\[P_{S} \df \{ x \in P\::\: \rank(x) \in S\} \cup \{ {\0},{\1}\}.\]
Denote by $f_S (P)$ the number of maximal chains of $P_S$. Equivalently,
$f_S(P)$ is the number of chains $x_1<\cdots<x_{|S|}$ in $P$ such that 
$\{\rank(x_1),\ldots,\rank(x_{|S|})\}=S$. The vector 
$\left(f_S (P)\::\: S\subseteq [1,n]\right) $
is called the {\em flag
$f$-vector} of $P$. Whenever it does not cause confusion, we write 
$f_{s_1\,\ldots\, s_k}$ rather than $f_{\{s_1,\ldots,s_k\}}$; in particular,
$f_{\{m\}}$ is always denoted $f_m$. 

Various properties of the flag $f$-vector are more easily seen in 
different bases. An often used equivalent encoding is the {\em flag 
$h$-vector} 
\mbox{$(h_S(P)\::\: S\subseteq [1,n])$} given by the formula
\[h_S (P)\df \sum _{T\subseteq S} (-1)^{|S\setminus T|} f_T(P),\]
or, equivalently,
\[f_S (P)= \sum _{T\subseteq S} h_T(P).\]
The {\em $ab$-index} $\Psi_P(a,b)$ of $P$ is a generating function
for the flag $h$-vector.
It is the following polynomial in the noncommuting variables $a$ and $b$:
\begin{equation}
\label{E_ab}
   \Psi_P(a,b)
  =
   \sum_{S \subseteq [1,n]} h_S (P)  u_S,
\end{equation}
where $u_S$ is the monomial $u_1 u_2\cdots u_n$ with $u_i=a$ if $i\not\in S$,
and $u_i=b$ if $i\in S$.

The {\em M\"obius function} of a graded poset                                   
$P$ is defined recursively for any subinterval of $P$ by the formula                                       
\[\mu ([x,y])=                                                                  
\left\{                                                                         
\begin{array}{cl}                                                               
1 & \mbox{if $x=y$},\\                                                          
-\sum_{x\le z<y} \mu([x,z]) &\mbox{otherwise}.\\   
\end{array}                                                                     
\right.                                                                         
\]
Equivalently, by Philip Hall's theorem,
the M\"obius function of a
graded poset $P$ of rank $n+1$ is the                                          
{\em reduced Euler characteristic} of the order complex, i.e.,  
it is given by the formula    
\begin{equation} 
\label{E_MEc}
\mu (P)=\sum_{S\subseteq [1,n]} (-1)^{|S|+1} f_S(P). 
\end{equation} 
(See \cite[Proposition 3.8.5]{Stanley-EC}.)                                     

A graded poset $P$ is {\em Eulerian} if the M\"obius function of every interval 
$[x,y]$ is given by $\mu([x,y]) = (-1)^{\rank(x,y)}$.                             
(Here $\rank(x,y)=\rank([x,y])=\rank(y) - \rank(x)$.)                        

The first characterization of all linear equalities holding 
for the flag $f$-vectors of all Eulerian posets was given by Bayer and 
Billera in \cite{Bayer-Billera}. 
The equations of the theorem are called the generalized Dehn-Sommerville 
equations.
Call the subspace of ${\bf R}^{2^n}$ they determine the {\em Eulerian subspace};
its dimension is the Fibonacci
number $e_n$ ($e_0=e_1=1$, $e_n=e_{n-1}+e_{n-2}$).

\begin{theorem}[Bayer and Billera]
\label{T_BB}
Every linear equality holding for the flag $f$-vector of all Eulerian posets
of rank $n+1$ is a consequence of the equalities
\[\left((-1)^{i-1}+(-1)^{k+1}\right)f_S+\sum _{j=i}^k (-1)^j
f_{S\cup\{j\}}=0\]
for $S\subseteq [1,n]$ and
$[i,k]$ a maximal interval of $[1,n]\setminus S$.
\end{theorem}

Fine discovered that the $ab$-index of a polytope can be written as a 
polynomial in the noncommuting variables $c\df a+b$ and $d\df ab + ba$.
Bayer and Klapper \cite{Bayer-Klapper} proved that for a graded poset $P$,
the equations of Theorem~\ref{T_BB} hold if and only if the $ab$-index is a 
polynomial with integer coefficients in $c$ and $d$.
This polynomial is called the {\em $cd$-index} of $P$.
Stanley (\cite{Stanley-flag})
gives an explicit recursion for the $cd$-index in terms of intervals 
of $P$ for Eulerian posets.
(He thus gives another proof of the existence of the $cd$-index for Eulerian
posets.)

\subsection{The flag \lv-vector and the flag  \LV-vector}
\label{s_ell}

The introduction of another vector equivalent to the flag $f$-vector
simplifies calculations.

\begin{definition}                                                              
\label{D_lv}                                                                    
{\em                                                                            
The {\em flag $\lv$-vector} of a graded partially ordered set $P$ of            
rank $n+1$ is the vector $(\lv_S(P) \:: \: S\subseteq [1,n])$,            
where                                                                           
\[\lv_S(P)\df                                                            
(-1)^{n-|S|} \sum_{T\supseteq [1,n]\setminus S}  
(-1)^{|T|} f_T(P).\]
}                                                                               
\end{definition}                                                                
As a consequence,                                                               
\begin{equation}                                                                
\label{E_fl}                                                                    
f_S(P)=\sum _{T\subseteq [1,n]\setminus S} \lv_T(P).               
\end{equation}                                                                  
The flag $\lv$-vector was first considered by Billera and Hetyei 
(\cite{Billera-Hetyei}) while          
describing all linear inequalities holding for the flag $f$-vectors
of all graded partially ordered sets.                                              
It turned out to give a sparse                             
representation of the cone of flag $f$-vectors described in that paper.

A variant significant for Eulerian posets is the flag $\LV$-vector.
\begin{definition}
\label{D_LV}
{\em
The {\em flag $\LV$-vector} of a graded 
partially ordered set $P$ of
rank $n+1$ is the vector $(\LV_S(P) \::\: S\subseteq [1,n])$, where
\[\LV_S (P)\df (-1)^{n-|S|}
\sum_{T\supseteq [1,n]\setminus S} \left(-{1\over 2}\right)^{|T|}
f_T(P).\]
}
\end{definition}

Inverting the relation of the definition gives
\[f_S(P)=2^{|S|}\sum _{T\subseteq [1,n]\setminus S} \LV_T(P).\]

When the poset $P$ is Eulerian,
the parameters $\LV_S (P)$ are actually the coefficients of the
{\em $ce$-index} of the poset $P$.
The $ce$-index was introduced by Stanley (\cite{Stanley-flag}) as an
alternative way of viewing the $cd$-index.
The letter $c$ continues to stand for $a+b$; now let $e\df a-b$.
The $ab$-index of a poset can be written in terms of $c$ and $d$ if and only
if it can be written in terms of $c$ and $ee$.
It is easy to verify that $\LV_S (P)$ is exactly the coefficient
in the \mbox{{\em $ce$-index}} of $P$ of the word 
$u_S=u_1u_2\cdots u_n$ where
$u_i=c$ if $i\not\in S$,
and $u_i=e$ if $i\in S$.
Since the existence of the $cd$-index is equivalent to the validity of the
generalized Dehn-Sommerville equations, we get the following proposition.
(It can be proved directly from the definition of the
flag $\LV$-vector, yielding an alternative way to prove the        
existence of the  $cd$-index for Eulerian posets.)
A subset $S\subseteq [1,n]$ is {\em even} if all the maximal intervals contained
in $S$ are of even length.

\begin{proposition}
\label{P_BBl}
The generalized Dehn-Sommerville relations hold for a poset $P$ if and only if 
$\LV_S(P)=0$ whenever $S$ is not an even set.
\end{proposition}  

The generalized Dehn-Sommerville relations hold (by chance) for some 
nonEulerian posets.  
A poset is Eulerian, however, if these relations hold for all intervals of 
the poset.

\begin{corollary} 
\label{P_EulerLV}                                                               
A graded partially ordered set is Eulerian if and only if 
$\LV_S ([x,y])=0$                                
for every             
interval $[x,y]\subseteq P$  and every subset $S$ of \mbox{$[1,\rank(x,y)-1]$}
that is not an even set.
\end{corollary}                                                               

\section{Half-Eulerian posets}
\label{half-Eul}

In this section we find special points in the closed cone of flag vectors
of Eulerian posets.
First consider the extremes of the closed cone of flag vectors of all graded
posets, found by Billera and Hetyei (\cite{Billera-Hetyei}).

\begin{definition}
{\em Given a graded 
poset $P$ of rank $n+1$, an interval $I\subseteq [1,n]$, and a positive  
integer $k$, 
$\HD_I^k(P)$ is the graded poset obtained from $P$ by replacing 
every $x\in P$ with rank in $I$ by $k$ elements $x_1,\ldots,x_k$
and by imposing the following relations.
\begin{enumerate}
\item[($i$)] If for  $x,y\in P$, $\rank(x)\in I$ and $\rank(y)\not\in I$,
then $x_i<y$ in $\HD_I^k(P)$ if and only if $x<y$ in $P$,
and $y<x_i$ in $\HD_I^k(P)$ if and only if $y<x$ in $P$.
\item[($ii$)] If $\{\rank (x),\rank (y)\}\subseteq I$, then $x_i<y_j$ 
in $\HD_I^k(P)$  if and only if $i=j$ and $x<y$ in $P$.
\end{enumerate}
}\end{definition}

Clearly $\HD_I^k P$  
is a graded poset of the same rank as $P$.
Its flag $f$-vector can be computed from that of $P$ in a straightforward
manner.

An {\em interval system on $[1,n]$} is any set of subintervals of $[1,n]$
that form an antichain (that is, no interval is contained in another).
(Much of what follows holds even if the intervals do not form an antichain,
but the assumption simplifies the statements of some theorems.)
For any interval system $\iv$ on $[1,n]$, and any positive integer $N$,
the poset $P(n,\iv,N)$ is defined to be the poset obtained from a chain
of rank $n+1$ by applying $\HD_I^N$ for all $I\in\iv$. 
It does not matter in which order these operators are applied. 
(Different values of $N$ can be used for each interval $I$, but we do not
need that generality  here.)
Consider the sequence of posets for a fixed interval system $\iv$ as 
$N$ goes to infinity.
Billera and Hetyei (\cite{Billera-Hetyei}) showed that the normalized flag 
vectors of such a sequence converge to a vector on an extreme ray of the cone of
flag vectors of all graded posets.
More precisely, 
\begin{theorem}[Billera and Hetyei]
\label{Billera-Hetyei theorem}
Suppose $\iv$ is an interval system of $k$ intervals on $[1,n]$.
Then the vector 
\[\left(\lim_{N\rightarrow \infty} {1\over N^k} 
f_S (P(n,\iv, N)):S\subseteq [1,n]\right)\]
generates an extreme ray of the cone of flag vectors of all graded posets.
Moreover, all extreme rays are generated in this way.
\end{theorem}

Unfortunately, none of the posets $P(n,\iv,N)$ are Eulerian, and 
none of these extreme rays are contained in the closed
cone of flag vectors of Eulerian posets.
However some of the posets are ``half-Eulerian'', and lead us to extreme
rays of the Eulerian cone.

For the interval system ${\cal I}=\{[1,1],[2,2],\ldots,[n,n]\}$,
abbreviate $\HD_{\cal I}^2(P)$ as $\HD P$, and call this the {\em horizontal 
double} of $P$.
Thus the horizontal double of $P$ is the poset obtained 
from $P$ by replacing every $x\in P\setminus \{\0,\1\}$
with two elements $x_1,x_2$ such that $\0$ and $\1$ remain the minimum 
and maximum elements of the partially ordered set, and 
$x_i<y_j$ if and only if $x<y$ in $P$. 
(In the Hasse diagram of $P$, every edge is replaced by $\Join$.)

\begin{definition}
{\em A {\em half-Eulerian poset} is a graded partially ordered set whose 
horizontal double is Eulerian.
}\end{definition}
For more information on half-Eulerian posets, see \cite{Bayer-Hetyei}.

The flag $f$-vectors of $P$ and its  
horizontal double are connected by the formula
$f_S (\HD P)=2^{|S|} f_S (P)$.
Thus, 
\begin{equation} \label{elltoL} 
\LV_S(\HD P)= \lv_S(P).\end{equation}

Applying the definition of Eulerian to the horizontal double of a poset
we get
\begin{proposition}
A graded partially ordered set $P$ is half-Eulerian if and only if for every
interval $[x,y]$ of $P$,
\[ \sum_{i=1}^{\rank(x,y)-1} (-1)^{i-1}  f_i ([x,y])=
 (1+(-1)^{\rank(x,y)})/2.\]
\end{proposition}

Corollary \ref{P_EulerLV}  can now be restated for half-Eulerian posets.
\begin{proposition}  
\label{D_hE}                                                                    
A graded partially ordered set is half-Eulerian if and only if 
$\lv_S ([x,y])=0$                                
for every             
interval $[x,y]\subseteq P$  and every subset $S$ of \mbox{$[1,\rank(x,y)-1]$}
that is not an even set.
\end{proposition}                                                               

The flag vectors of the horizontal doubles of half-Eulerian posets span
the Eulerian subspace, the subspace defined by the generalized 
Dehn-Sommerville equations.
But the cones they determine may be different.
Recall ${\cal C}^{n+1}_{\cal E}$ is the closed cone of 
flag vectors of Eulerian posets.
Now write ${\cal C}^{n+1}_{\cal D}$ for the closed cone of 
flag vectors of horizontal doubles of half-Eulerian posets.
We do not know if the inclusion 
${\cal C}^{n+1}_{\cal D}\subseteq {\cal C}^{n+1}_{\cal E}$ is actually 
equality.

For which interval systems $\iv$ is $P(n,\iv,N)$ half-Eulerian?

\begin{definition}                                                     
{\em                                                              
An interval system $\iv$ on $[1,n]$ is {\em even} if for every
pair of intervals $I,J\in \iv$ the intersection $I\cap J$ has
an even number of elements. (In particular, $|I|$ must be even for every
$I\in \iv$.) 
}                                  
\end{definition}

Our goal is to show that the posets $P(n,\iv,N)$ are half-Eulerian if and 
only if $\iv$ is an even interval system.
For this we need to understand the intervals of the posets $P(n,\iv,N)$.

\begin{proposition}
\label{P_limitint}
The interval $[x,y]\subseteq P(n,\iv,N)$ is isomorphic to \\
$P(\rank(x,y)-1,{\cal J},N)$, where 
${\cal J}=\{I-\rank(x)\::\: 
I\in \iv,\, I\subseteq [\rank (x)+1, \rank(y)-1] \} 
$.
\end{proposition}
\begin{proof}
Let $\rank(x)=r$ and $\rank (y)=s$. 
Construct $P(n,\iv,N)$  
by applying the operators $\HD_I^N$ for all $I\in\iv$ to a chain.
Since the order of applying these operators is arbitrary, we may
choose to apply first those for which $I$ is not a subset of 
$[r+1,s-1]$. 
At this point for every $x'$ of rank $r$ and $y'$ of rank $s$ with  $y'\ge x'$,
the interval $[x',y']$ is isomorphic to a chain of rank 
$\rank (x',y')$. Applying the remaining operators $\HD_I^N$ leaves 
the elements of rank at most $r$ or of rank at least $s$ unchanged, and 
has the same effect on $[x',y']$ as applying the operators $\HD_{I-r}^N$ 
to a chain of rank $\rank (x',y')$.  
\end{proof}

The effect on the flag $f$-vector of applying the operator $\HD_I^N$
to a poset of rank $n+1$ is given by the formula 
\begin{equation}
\label{E_HDINf}
f_S (\HD_I^N(P))
=\left\{
\begin{array}{ll}
N f_S(P) & \mbox{if $I\cap S\neq \emptyset$,}\\
f_S(P) &\mbox{otherwise.}\\
\end{array}
\right.
\end{equation}
This enables us to write an \lv-vector formula.

\begin{lemma}
\label{E_HDINl}
For $P$ a graded poset of rank $n+1$, $S\subseteq [1,n]$, and $N$ a positive
integer,
\begin{equation}
\lv_S (\HD_I^N(P))
=N \lv_S(P)-(N-1)\sum_{T\cup I=S} \lv_T(P). 
\end{equation}
\end{lemma}
\begin{proof} From the definition of $\ell_S$ and equation~(\ref{E_HDINf}),
\begin{eqnarray*}
\lv_S(\HD_I^N(P))
&=&(-1)^{n-|S|} \sum_{R\supseteq [1,n]\setminus S} 
(-1)^{|R|} f_R (\HD_I^N(P))\\  
&=&(-1)^{n-|S|} \sum_{R\supseteq [1,n]\setminus S} 
\!\!\! (-1)^{|R|} N f_R (P)  \\
& &
 \mbox{} -(-1)^{n-|S|} \sum_{R\supseteq [1,n]\setminus S \atop 
   R\subseteq [1,n]\setminus I} 
\!\!\!(-1)^{|R|} (N-1) f_R (P)\\
&=& N \lv_S(P)
 -(-1)^{n-|S|} \sum_{R\supseteq [1,n]\setminus S \atop 
   R\subseteq [1,n]\setminus I} 
(-1)^{|R|} (N-1) f_R (P)\\
\end{eqnarray*}
By (\ref{E_fl}), the coefficient in 
$-(-1)^{n-|S|} \sum_{R\supseteq [1,n]\setminus S \atop 
   R\subseteq [1,n]\setminus I} (-1)^{|R|} (N-1) f_R (P)$ 
of $\lv_T(P)$ is
\[ -(N-1)(-1)^{n-|S|} \sum_{R\supseteq [1,n]\setminus S \atop 
   R\subseteq [1,n]\setminus (T\cup I)} 
(-1)^{|R|}\, , \]
which is an empty sum if $(T\cup I)$ is not contained in 
$S$, zero if $(T\cup I)$
is properly contained in $S$, and 
$-(N-1)(-1)^{n-|S|}(-1)^{|[1,n]\setminus S|}=-(N-1)$ if 
$(T\cup I)=S$. 
This gives the recursion of the lemma.
\end{proof}

 From this we can determine which of the posets $P(n,\iv, N)$ are half-Eulerian. 

\begin{proposition}
\label{P_limitHE}
Let $\iv$ be an interval system on $[1,n]$.
\begin{enumerate}
\item If $\iv$ is an even system of intervals, then for all $N$ the partially
      ordered set $P(n,\iv,N)$ is half-Eulerian.
\item If for some $N>1$, $P(n,\iv,N)$ is half-Eulerian, then
      $\iv$ is an even system of intervals.
\end{enumerate}
\end{proposition}
\begin{proof}
Using Lemma~\ref{E_HDINl} we
can show by induction on $|\iv|$ that for every $N$,
$\lv^{n+1}_S\left(P(n,\iv, N)\right)$ 
is zero unless $S$ is the union of some intervals of $\iv$.
In particular, if $\iv$ is an even system of intervals, then 
$\lv_S\left(P(n,\iv, N)\right)=0$ whenever $S$ is not an even set.
The same observation holds for every
interval $[x,y]\subseteq P(n,\iv, N)$ as well, since by 
Proposition \ref{P_limitint} $[x,y]$ is isomorphic to $P(m,{\cal J}, N)$ for 
some  
$m\leq n$ and some even system of intervals $\cal J$. Therefore the 
conditions of Proposition~\ref{D_hE} are satisfied by $P(n,\iv, N)$ for
every $N$, if $\iv$ is an even system of intervals.

Now assume $\iv$ is a system of intervals that is not even.
First consider the case where $\iv$  
contains an interval $I_m=[a,b]$ with $b-a$ even (hence $I_m$
is odd).
Let ${\cal J}=\{I_m-a+1\}=\{[1,b-a+1]\}$.
For $S$ nonempty, $f_S(P(b-a+1,{\cal J},N))=N$, so
\begin{eqnarray*}
\lefteqn{\lv_{[1,b-a+1]}(P(b-a+1,{\cal J},N))}\hspace*{.5in}\\
&=& \sum_{T\subseteq[1,b-a+1]} (-1)^{|T|}f_T(P(b-a+1,{\cal J},N))\\
   &=& 1+\sum_{{T\subseteq[1,b-a+1]}\atop{T\ne\emptyset}} (-1)^{|T|}N=1-N.
\end{eqnarray*}
So $\lv_{[1,b-a+1]}(P(b-a+1,{\cal J},N))\ne 0$ for $N>1$.
Fix $N>1$, and choose $x$ and $y$ in $P(n,\iv,N)$ with $\rank(x)=a-1$,
$\rank(y)=b+1$, and $x\le y$.
Then by Proposition~\ref{P_limitint}, 
$\lv_{[1,\rank(x,y)-1]}([x,y])
=\lv_{[1,b-a+1]}(P(b-a+1,{\cal J},N))\ne 0$, with $|[1,b-a+1]|$ odd.
So $P(n,\iv,N)$ is not half-Eulerian.

Now suppose $\iv$ contains only even intervals, but some two intervals have an 
odd overlap. 
Let $I_p=[a,d]$ and $I_q=[c,b]$, where $a<c\le d<b$ and $d-a$ and $b-c$ are odd,
but $d-c$ is even.                                                              
Then $b-a$ is also even.                                                        
We show that we may assume no other interval of $\iv$ is in the union 
$I_p\cup I_q$.
Suppose $I_r=[e,f]$ is another interval of $\iv$ with $[e,f]\subset [a,b]$
(and $f-e$ is odd).
Since $\iv$ is an antichain, $a<e<c\le d<f<b$.
If $e-a$ is even, then 
$|I_q\cap I_r|=|[c,f]|=f-c+1 = (f-e)+(e-a)-(d-a)+(d-c)+1$,
which is odd, because it is the sum of three odds and two evens.
If $e-a$ is odd, then $|I_p\cap I_r|=|[e,d]|=d-e+1=(d-a)-(e-a)+1$, which is odd 
because it is the sum of three odds.
Thus, if two intervals of $\iv$ have odd intersection and their union contains a
third interval of $\iv$, then two intervals of $\iv$ with smaller union have 
odd intersection.

So we may assume $I_p=[a,d]$ and $I_q=[c,b]$ have odd intersection, and their 
union $[a,b]$ contains no other interval of $\iv$.
Let ${\cal J}
=\{I_p-a+1, \linebreak I_q-a+1\}
=\{[1,d-a+1],[c-a+1,b-a+1]\}
$.                                            
Then                                                                            
\begin{eqnarray*}
\lefteqn{
f_S(P(b-a+1,{\cal J},N))}\hspace*{.5in}\\
&=&\left\{\begin{array}{ll}                     
             1 & \mbox{if $S=\emptyset$}\\                                      
	     N^2 & \mbox{if $S\cap(I_p-a+1)\ne\emptyset$ and 
			    $S\cap(I_q-a+1)\ne\emptyset$}\\
             N & \mbox{otherwise.}                                              
             \end{array} \right. 
\end{eqnarray*}
So                                                                              
\begin{eqnarray*}                                                               
\lefteqn{                                                                       
\lv_{[1,b-a+1]}(P(b-a+1,{\cal J},N))}\hspace*{.25in}\\
&=&                                   
   \sum_{T\subseteq[1,b-a+1]} (-1)^{|T|}f_T(P(b-a+1,{\cal J},N))\\     
   &=& \sum_{T\subseteq[1,b-a+1]} (-1)^{|T|}N^2                                 
   + \sum_{T\subseteq[1,c-a]} (-1)^{|T|}(N-N^2)\\                               
   &+& \sum_{T\subseteq[d-a+2,b-a+1]} (-1)^{|T|}(N-N^2)                         
   + (1-2N+N^2) = (1-N)^2.                                                      
\end{eqnarray*}                        
So $\lv_{[1,b-a+1]}(P(b-a+1,{\cal J},N))\ne 0$ for $N>1$.               
Fix $N>1$, and choose $x$ and $y$ in $P(n,\iv,N)$ with $\rank(x)=a-1$,
$\rank(y)=b+1$, and $x\le y$.
Then by Proposition~\ref{P_limitint}, 
$\lv_{[1,\rank(x,y)-1]}([x,y])
=\lv_{[1,b-a+1]}(P(b-a+1,{\cal J},N))\ne 0$, with $|[1,b-a+1]|$ odd.
So $P(n,\iv,N)$ is not half-Eulerian.
\end{proof}

As will be seen later, even interval systems give rise to extreme rays of the 
cone of flag vectors of Eulerian posets.  
It is of interest, therefore, to count them.

\begin{proposition}
The number of even interval systems on $[1,n]$ is 
${n \choose \lfloor n/2\rfloor}$.
\end{proposition}
\begin{proof}
We define a one-to-one correspondence between even interval systems on $[1,n]$
and sequences $\lambda=(\lambda_1,\lambda_2,\ldots, \lambda_n)\in \{-1,1\}^n$
satisfying $\sum_i \lambda_i = 0$ if $n$ is even and 
$\sum_i \lambda_i = 1$ if $n$ is odd.
Clearly there are ${n \choose \lfloor n/2\rfloor}$ such sequences.

For $\iv$ an even interval system, define 
$\lambda(\iv)=(\lambda_1,\lambda_2,\ldots, \lambda_n)\in \{-1,1\}^n$,
where $\lambda_i = (-1)^i$ if $i$ is an endpoint of an interval
of $\iv$, and $\lambda_i = (-1)^{i-1}$ otherwise.
(Note that for an even interval system, no number can be an endpoint of
more than one interval.)
For \iv\ an even interval system, summing $(-1)^i$ over the endpoints of
intervals gives 0.  So
\begin{eqnarray*}
\sum_{i=1}^n \lambda_i &=& \sum_{i=1}^n (-1)^{i-1} +  
			   \sum_{{\mbox{{\scriptsize $i$ endpoint}}}\atop
			     {\mbox{{\scriptsize of interval}}}} 
                           2(-1)^i\\
                       &=& 
			     \sum_{i=1}^n(-1)^{i-1} 
		       = \left\{\begin{array}{rl}
				0 & \mbox{if $n$ is even}\\
				1 & \mbox{if $n$ is odd} \end{array}\right. .
\end{eqnarray*}

On the other hand, given a 
sequence $\lambda=(\lambda_1,\lambda_2,\ldots, \lambda_n)\in \{-1,1\}^n$
satisfying $\sum_i \lambda_i = 0$ if $n$ is even and 
$\sum_i \lambda_i = 1$ if $n$ is odd, construct an even interval system as
follows.
Let $s_1 < s_2 < \cdots < s_k$ be the sequence of indices $s$ for which
$\lambda_s = (-1)^s$.
Then 
$ \sum_{i=1}^n (-1)^{i-1} =\sum_{i=1}^n \lambda_i
=\sum_{i=1}^n (-1)^{i-1} + \sum_{j=1}^k 2(-1)^{s_j}$, so
$ \sum_{j=1}^k (-1)^{s_j}=0$.
Thus the sequence of $s_j$'s contains the same number of even numbers as odd.
Construct an interval system $\iv = \{[a_1,b_1], [a_2,b_2], \ldots, [a_m,b_m]\}$
($2m=k$) recursively as follows.
Let $a_1=s_1$ and let $b_1=s_j$ where $j$ is the least index such that $s_1$
and $s_j$ are of opposite parity.
Then $\iv = {[a_1,b_1]}\cup \iv'$, where $\iv'$ is the interval system 
associated with $s_2 < s_3 < s_4 <  \cdots < s_k$ with $b_1=s_j$ removed.
Clearly $[a_1,b_1]$ is of even length.
If $[a_1,b_1]\cap [a_i,b_i]\ne \emptyset$ for some interval $[a_i,b_i]$ of 
$\iv'$, then $a_i<b_1$, so by the choice of $b_1$, $a_i$ has the same parity
as $a_1$.  Thus $[a_1,b_1]\cap [a_i,b_i]=[a_i,b_1]$ is of even length.
Furthermore, $b_i$ and $b_1$ are of the same parity, since $a_i$ and $a_1$ are,
so again by the choice of $b_1$, $b_i>b_1$.  
So the interval $[a_i,b_i]$ is not contained in the interval $[a_1,b_1]$.
The interval system $\{[a_m,b_m]\}$, is even, so by induction 
$\iv$ is an even interval system.

These constructions are inverses, giving the desired bijection.
\end{proof}

Recall that Billera and Hetyei (\cite{Billera-Hetyei}) found extremes of the 
cone of flag vectors of graded
posets as limits of the normalized flag vectors of the posets $P(n,\iv,N)$.
The next proposition follows easily by induction from Lemma~\ref{E_HDINl}.
\begin{proposition}
\label{P_limitl}
Let $\iv=\{I_1,I_2,\ldots,I_k\}$ be a system of $k\geq 0$ intervals on~$[1,n]$.
Then
\begin{eqnarray*}
\lefteqn{
\lim_{N\longrightarrow \infty} {1\over N^k} \lv_S 
\left(P(n,\iv, N)\right)} \hspace*{.5in}\\
&=&\sum _{j=0}^k (-1)^j 
\left|\left\{1\leq i_1<\cdots <i_j\leq k \::\: 
I_{i_1}\cup \cdots \cup I_{i_j}=S
\right\}\right|. 
\end{eqnarray*}
\end{proposition}

Write $f_S(P(n,\iv))=
\lim_{N\rightarrow\infty} f_S(P(n,\iv,N))/N^{|\iv|}$.
The vector these form (as $S$ ranges over all subsets of $[1,n]$) is not the
flag $f$-vector of an actual poset, but it is in the closed cone of flag
$f$-vectors of all graded posets.
We call the symbol $P(n,\iv)$ a ``limit poset'' and refer to the flag vector
of the limit poset.
If $\iv$ is an even interval system, then 
$(f_S(P(n,\iv)): S\subseteq [1,n])$ is in the closed cone of flag vectors
of half-Eulerian posets.
To get Eulerian posets the horizontal double operator is applied to 
$P(n,\iv,N)$.
The vector $(f_S(\HD P(n,\iv)): S\subseteq [1,n])$ is defined as a limit 
of the resulting normalized flag $f$-vectors, and satisfies
$f_S(\HD P(n,\iv))=2^{|S|}f_S(P(n,\iv))$.
It lies in the cone ${\cal C}^{n+1}_{{\cal D}}$ of flag vectors of doubles of 
half-Eulerian posets, a subcone of the Eulerian cone.

Recall (equation~(\ref{elltoL})) that the \lv-vector of a poset $P$ equals the
\LV-vector of its horizontal double $\HD P$.
The same holds after passing to the limit posets.
Thus, Proposition~\ref{P_limitl} gives 
\[ L_S(\HD P(n,{\cal I}))
=\sum _{j=0}^k (-1)^j 
\left|\left\{1\leq i_1<\cdots <i_j\leq k :\: 
I_{i_1}\cup \cdots \cup I_{i_j}=S
\right\}\right|, \]
where $\iv=\{I_1,I_2,\ldots,I_k\}$.

We look at the associated $cd$-indices of the ``doubled limit posets.''
Think of a word in $c$ and $d$ as a string with each $c$ occupying one position
and each $d$ occupying two positions.
The {\em weight} of a $cd$-word $w$ is then the number of positions of the 
string.
Associated to each $cd$-word $w$ is the even set $S(w)$ consisting of the 
positions occupied by the $d$'s.

\begin{proposition}
\label{cd}
For each $cd$-word $w$ with $k$ $d$'s and weight $n$, there exists an even 
interval system $\iv_w$ for which the $cd$-index of $\HD P(n,\iv_w)$ is $2^k w$.
\end{proposition}
\begin{proof}
Fix a $cd$-word $w$ with $k$ $d$'s and weight $n$.
Write the elements of $S(w)$ in increasing order as $i_1$, $i_1+1$, $i_2$,
$i_2+1$, \ldots, $i_k$, $i_k+1$, and let $\iv_w$ be the interval system
$\{[i_1, i_1+1], [i_2, i_2+1], \ldots, [i_k, i_k+1]\}$.
Let $\Phi=2^k w$.  Rewrite the $cd$-polynomial $\Phi$ as a $ce$-polynomial.
Recall from Sections~\ref{ss_gpo} and~\ref{s_ell} that $c=a+b$, $d=ab+ba$, and
$e=a-b$, so $d = (cc-ee)/2$.
Thus, $\Phi$ is rewritten as a sum of $2^k$ terms. Each is the result of
replacing some subset of the $d$'s by $cc$, and the rest by $ee$; the
coefficient is $\pm 1$, depending on whether the number of $d$'s
replaced by $ee$ is even or odd.
Thus
\[ 2^k w = \sum_{J\subseteq [1,k]} (-1)^J w_J,\]
where $w_J=w_1w_2\cdots w_n$, with $w_{i_j}=w_{i_j+1}=e$ if $j\in J$
and the remaining $w_i$'s are $c$.
By the \LV-vector version of Proposition~\ref{P_limitl}, this is precisely the 
$ce$-index of $\HD P(n,\iv_w)$.
\end{proof}

In \cite{Stanley-flag} Stanley first found for each $cd$-word $w$ a
sequence of Eulerian posets whose normalized $cd$-indices converge to $w$.
Our limit posets are closely related to Stanley's, but this particular
construction highlights the important link between the half-Eulerian and
Eulerian cones.

Before turning to inequalities satisfied by the flag vectors of Eulerian
posets, we consider the question of whether the two cones 
${\cal C}^{n+1}_{\cal D}$ and ${\cal C}^{n+1}_{\cal E}$ are equal.
For low ranks the two cones are the same, as seen below.
We know of no example in any rank of an Eulerian poset whose flag vector
is not contained in the cone ${\cal C}^{n+1}_{\cal D}$
of doubled half-Eulerian posets.
To look for such an example we turn to 
the best known examples of Eulerian posets, the face lattices of polytopes.
In \cite{Stanley-flag} Stanley proved the nonnegativity of the
$cd$-index for ``$S$-shellable regular CW-spheres'', a class of Eulerian
posets that includes all polytopes.  
By a result of Billera, Ehrenborg, and Readdy (\cite{Billera-Ehrenborg-Readdy}),
the lattice of regions of any oriented matroid also has a nonnegative 
$cd$-index.
Proposition~\ref{cd} implies that nonnegative
$cd$-indices (and the associated flag vectors) are in the 
cone generated by the $cd$-indices (flag vectors) of the doubles of limit posets
associated with even interval systems.
\begin{corollary}
${\cal C}^{n+1}_{\cal D}$ 
contains the flag vectors of all Eulerian  posets 
with nonnegative $cd$-indices.
This includes the face lattices of polytopes and the lattices of regions of
oriented matroids.
\end{corollary}

\begin{conjecture}
The closed cone ${\cal C}^{n+1}_{\cal E}$ of flag vectors of Eulerian posets is
the same as the closed cone ${\cal C}^{n+1}_{\cal D}$
of flag vectors of horizontal doubles of half-Eulerian posets.
\end{conjecture}

\section{Inequalities}
\label{sec-inequalities}

Throughout this section we use the following notation.
\begin{definition}
{\em
The {\em interval system $\iv[S]$ of a set $S\subseteq [1,n]$}  is the
family of intervals  
$\iv[S]=\{[a_1,b_1],\ldots,[a_k,b_k]\}$, where
$S=[a_1,b_1]\cup \cdots\cup [a_k,b_k]$ and
$b_{i-1}<a_i-1$  for $i\geq 2$. 
In other words, $\iv[S]$ is the collection of the maximal intervals contained
in $S$.   
}
\end{definition}
Note that $S$ is an even set if and only if $\iv[S]$ is an even interval system.

The following flag vector forms can be proved nonnegative by writing them as
convolutions of basic nonnegative forms \cite{Billera-Liu,Kalai}. (See
\ref{S_ring}.)  
The issue of whether they give all linear inequalities on flag vectors of 
Eulerian posets was raised by Billera and Liu (see the discussion after
Proposition~1.3 in \cite{Billera-Liu}).
We give here a simple direct argument for their nonnegativity that avoids
convolutions. 

\begin{proposition}[Inequality Lemma]
\label{ineqlemma}
Let $T$ and $V$ be subsets of $[1,n]$ such that for every 
$I\in \iv[V]$, $|I\cap T|\le 1$.
Write $S=[1,n]\setminus V$.
For $P$ any rank $n+1$ Eulerian poset,
\[\sum_{R\subseteq T}(-2)^{|T\setminus R| } f_{S\cup R}(P)\ge 0. \]
Equivalently,
\[ (-1)^{|T|}\sum_{T\subseteq Q\subseteq V} L_Q(P)\ge 0.\] 
\end{proposition}
\begin{proof}
The idea is that since no two elements of $T$ are in the same gap of $S$,
elements with ranks in $T$ can be inserted independently in chains with rank 
set $S$.
For $C$ an $S$-chain (i.e., a chain with rank set $S$) and $t\in T$, let 
$n_t(C)$ be the number of rank $t$ elements $x\in P$ such that $C\cup \{x\}$ 
is a chain of $P$.
Since every interval of an Eulerian poset is Eulerian, $n_t(C)\ge 2$ for all
$C$ and $t$.
So
\begin{eqnarray*}
\sum_{R\subseteq T} (-2)^{|T\setminus R|} f_{S\cup R} (P)
&=& \sum_{R\subseteq T} (-2)^{|T\setminus R|} 
    \sum_{\mbox{\scriptsize $C$ an $S$-chain}} \prod_{t\in R} n_t(C)\\
&=& \sum_{\mbox{\scriptsize $C$ an $S$-chain}} 
    \sum_{R\subseteq T} (-2)^{|T\setminus R|} \prod_{t\in R} n_t(C)\\
&=& \sum_{\mbox{\scriptsize $C$ an $S$-chain}}  
    \prod_{t\in T} (n_t(C)-2)\ge 0.\\
\end{eqnarray*}
So the flag vector inequality is proved.
The second inequality is simply the translation into \LV-vector form.
\end{proof} 

Here are some new inequalities.

\begin{theorem}
\label{ijk-ineqs}
Let $1 \le i < j < k \le n$.
For $P$ any rank $n+1$ Eulerian poset, 
\[ f_{ik}(P) - 2f_i(P) - 2f_k(P) + 2f_j(P) \ge 0. \]
\end{theorem}
\begin{proof}
First order the rank $j$ elements of $P$ in the following way.
Choose any order, $G_1$, $G_2$, \ldots, $G_m$ for the components of the Hasse
diagram of the rank-selected poset $P_{\{i,j,k\}}$.
For each rank $j$ element $y$ of $P$, identify the component containing $y$ 
by $y\in G_{g(y)}$.
Order the rank $j$ elements of $P$ in any way consistent with the ordering of
components.
That is, choose an order $y_1$, $y_2$, \ldots, $y_r$ such that $y_s<y_t$
implies $g(y_s)\le g(y_t)$.

A rank $i$ element $x$ {\em belongs} to $y_q$ if $q$ is the least index such
that $x < y_q$ in $P$.
Write $I_q$ for the number of rank $i$ elements belonging to $y_q$, and
$I'_q$ for the number of rank $i$ elements $x$ such that $x < y_q$, but $x$
does not belong to $y_q$.
Similarly, 
a rank $k$ element $z$ {\em belongs} to $y_q$ if $q$ is the least index such
that $y_q  < z$ in $P$.
Write $K_q$ for the number of rank $k$ elements belonging to $y_q$, and
$K'_q$ for the number of rank $k$ elements $z$ such that $y_q  < z$, but $z$
does not belong to $y_q$.
Note that $I_q + I'_q \ge 2$ and $K_q + K'_q \ge 2$, since $P$ is Eulerian.
A flag $x<z$ {\em belongs} to $y_q$ if $x<y_q<z$ and $q$ is the least index such
that either $x<y_q$ or $y_q<z$.

Let $F= f_{ik}(P) - 2f_i(P) - 2f_k(P) + 2f_j(P)$.
Let $F_q$ be the contribution to $F$ by elements and flags belonging to $y_q$.
Thus, 
\[F_q = I_q K_q + I'_q K_q + I_q K'_q - 2I_q - 2K_q + 2.\]

If $I'_q\ge 2$, then $F_q=I_q(K_q+K'_q-2)+(I'_q-2)K_q+2 \ge 2$.

If $I'_q=K'_q=0$, then $F_q=(I_q-2)(K_q-2)-2\ge -2$.

In all other cases it is easy to check that $F_q\ge 0$.

Suppose that the rank $j$ elements in component $G_\ell$ are 
$y_s$, $y_{s+1}$, \ldots, $y_t$.
Then $I'_s = K'_s = 0$, so $F_s\ge -2$.
Furthermore, $I_t = K_t = 0$, because any rank $i$ element $x$ related to 
$y_t$ must also be related to at least one other rank $j$ element, and it is in 
the same component.
That rank $j$ element has index less than $t$, so $x$ does not belong to $y_t$.
This in turn implies $I'_t\ge 2$, so $F_t\ge 2$.
For all $q$, $s < q < t$, either $I'_q>0$ or $K'_q>0$, by the connectivity
of the component, so $F_q\ge 0$.
Thus $\sum_{q=s}^t F_q \ge 0$.
This is true for each component $G_\ell$, so $F=\sum_{q=1}^r F_q \ge 0$.
\end{proof}

These inequalities can be used to generate others by convolution (see
~\ref{S_ring}.)  

Evaluating the flag vector inequalities of Proposition~\ref{ineqlemma} for the
horizontal double $\HD P$ of a half-Eulerian poset $P$ gives the inequalities,
for $S$ and $T$ satisfying the hypotheses of Proposition~\ref{ineqlemma},
\begin{equation}
\label{halfEulineq1}
\sum_{R\subseteq T}(-1)^{|T\setminus R| } f_{S\cup R}(P)\ge 0. 
\end{equation}
These inequalities are valid not just for half-Eulerian posets but for {\em all}
graded posets.
The proof of Proposition~\ref{ineqlemma} uses only the fact that in every open
interval of an Eulerian poset there are at least two elements of each rank.
If the proof is rewritten using the assumption that in every open interval
there is at least one element of each rank, the 
inequalities~(\ref{halfEulineq1}) are proved for all  graded posets.

Similarly, the flag vector inequalities of Theorem~\ref{ijk-ineqs} give
inequalities for half-Eulerian posets,
\[ f_{ik}(P) - f_i(P) - f_k(P) + f_j(P) \ge 0. \]
The proof of Theorem~\ref{ijk-ineqs} can be modified in the same way to show
these inequalities are valid for all graded posets.
The first instance of this class of inequalities was found by 
Billera and Liu (\cite{Billera-Liu}).

We conjecture that all inequalities valid for half-Eulerian posets come from
inequalities valid for all graded posets.
Inequalities for half-Eulerian posets are to be interpreted as conditions in
the subspace of ${\bf R}^{2^n}$ spanned by flag vectors of half-Eulerian
posets, but we are describing them in ${\bf R}^{2^n}$.
Giving inequalities using linear forms in the flag numbers $f_S$ over 
${\bf R}^{2^n}$, the statement is as follows.
\begin{conjecture}
Every linear form that is nonnegative for the flag vectors of all half-Eulerian
posets is the sum of a linear form that is nonnegative for all graded posets and
a linear form that is zero for all half-Eulerian posets.
\end{conjecture}

\section{Extreme Rays and Facets of the Cone}
\label{sec-cone}

We have described some points in the Eulerian cone ${\cal C}^{n+1}_{{\cal E}}$
and some inequalities satisfied by all points in the cone.
We turn now to identifying which of these give extreme rays and facets.

If $\iv$ is an even interval system, then 
$(f_S(P(n,\iv)): S\subseteq [1,n])$ 
is on an extreme ray in the closed
cone of flag vectors of all graded posets, and is in the subcone of flag 
$f$-vectors of half-Eulerian posets. 
Therefore it is on an extreme ray of the subcone.

\begin{proposition}
For every even interval system $\iv$, the flag vector of the limit poset
$P(n,\iv)$ generates an extreme ray of the cone of flag vectors of half-Eulerian
posets.
\end{proposition}

What does this say about the extreme rays of the cone of flag vectors of
Eulerian posets?
For every even interval system \iv,
the flag vector of $\HD P(n,\iv)$ lies on an extreme ray of the subcone
${\cal C}^{n+1}_{{\cal D}}$, but we
cannot conclude directly that it lies on an extreme ray
of the cone ${\cal C}^{n+1}_{{\cal E}}$.
A separate proof is needed.

For the following proofs, we use the computation of $\lv_Q(P(n,\iv))$
(and $\LV_Q(\HD P(n,\iv))$)
from the decompositions of $Q$ as the union of intervals of $\cal I$
(Proposition~\ref{P_limitl}).

\begin{theorem}
For every even interval system $\iv$, the flag vector of the doubled limit poset
$\HD P(n,\iv)$ generates an extreme ray of the cone of flag vectors of 
Eulerian posets.
\end{theorem}
\begin{proof}
We work in the closed cone of \LV-vectors of Eulerian posets.
The cone of \LV-vectors of Eulerian posets is contained in the subspace of 
${\bf R}^{2^n}$ determined by the equations $L_S=0$ for $S$ not an even set.
To prove that  the \LV-vector of $\HD P(n,\iv)$ generates an extreme ray,
we show that it lies on linearly independent supporting hyperplanes,
one for each nonempty even set $V$ in $[1,n]$.
Fix an even interval system $\iv$.
For each nonempty even set $V\subseteq[1,n]$, we find a set $T$ such that
$T$ and $V$ satisfy the hypothesis of Proposition~\ref{ineqlemma} and 
$\sum_{T\subseteq Q\subseteq V} L_Q(\HD P(n,\iv))=0$.

{\em Case 1.} Suppose $V$ is the union of some intervals in \iv. 
Let $I_1$, $I_2$, \ldots, $I_k$ be all the intervals of \iv\ contained in $V$.
Set $T=\emptyset$.
Then for each subset $J\subseteq [1,k]$, the corresponding union of intervals 
contributes $(-1)^{|J|}$ to $L_Q(\HD P(n,\iv))$, for $Q=\cup_{j\in J} I_j$.
Thus 
$\sum_{T\subseteq Q\subseteq V} L_Q(\HD P(n,\iv))=
\sum_{J\subseteq [1,k]}(-1)^{|J|}=0$.

{\em Case 2.} If $V$ is not the union of some intervals in \iv, let $W$ be the
union of all those intervals of \iv\ contained in $V$.
Choose $t\in V\setminus W$, and set $T=\{t\}$.
For $Q\subseteq V$, $L_Q(\HD P(n,\iv))=0$ unless $Q\subseteq W$. 
But if $Q\subseteq W$ then $t$ cannot be in $Q$.
So 
$\sum_{\{t\}\subseteq Q\subseteq V} L_Q(\HD P(n,\iv))=0$.

Now $\sum_{T\subseteq Q\subseteq V} L_Q(P)=0$ determines a 
supporting hyperplane of the closed cone of \LV-vectors of Eulerian posets,
because the inequality of Proposition~\ref{ineqlemma} is valid, and the
poset $\HD P(n,\iv)$ lies on the hyperplane.
The hyperplane equations each involve a distinct maximal set $V$, which is even,
so they are linearly independent on the subspace determined by the equations
$L_S=0$ for $S$ not an even set.
So the doubled limit poset $\HD P(n,\iv)$ is on an extreme ray of the cone.
\end{proof}

Note how far we are, however, from a complete description of the extreme
rays.

\begin{conjecture}
For every positive integer $n$, the closed cone of flag $f$-vectors of 
Eulerian posets of rank $n+1$ is finitely generated.
\end{conjecture}

\begin{lemma}[Facet Lemma]
\label{facetlemma}
Assume $\sum_{Q\subseteq[1,n]} a_QL_Q(P)\ge 0$ for all 
Eulerian posets $P$ of 
rank $n+1$.
Let $M\subseteq [1,n]$ be a fixed even set.
Suppose for all even sets $R\subseteq [1,n]$, $R\ne M$, there exists an
interval system $\iv(R)$ consisting of disjoint even intervals whose union is
$R$ and such that 
$\sum_{Q\subseteq[1,n]} a_Q\lv_Q(P(n,\iv(R)))= 0$.
Then 
$\sum_{Q\subseteq[1,n]} a_QL_Q(P)= 0$ determines a facet of the
closed cone of \LV-vectors of Eulerian posets.
\end{lemma}
(Note that $\iv(R)$ need not be $\iv[R]$.)

\begin{proof}
The dimension of the cone ${\cal C}^{n+1}_{\cal E}$ equals the number of even 
subsets (a Fibonacci number).
So it suffices to show that the vectors 
$(\lv_Q(P(n,\iv(R))))$ $\mbox{}=
(L_Q(\HD P(n,\iv(R) )))$
are linearly independent.
To see this, note that
for every set $Q$ not contained in $R$, $\lv_Q(P(n,\iv(R)))=0$.
By the disjointness of the intervals in $\iv(R)$, there is a unique way to
write $R$ as the union of intervals in $\iv(R)$.
So by Proposition~\ref{P_limitl}, $(\lv_R(P(n,\iv(R))))=(-1)^{|\iv(R)|}$.
Thus, $R$ is the unique maximal set $Q$ for which 
$(\lv_Q(P(n,\iv(R))))\ne 0$.
So the \LV-vectors of the posets $\HD P(n,\iv(R))$, as $R$ ranges over sets
different from $M$, are linearly independent.
\end{proof}

\begin{proposition}
\label{ineq0}
The inequality 
$\sum_{Q\subseteq[1,n]} L_Q(P)\ge 0$ (or, equivalently, \linebreak
$f_\emptyset(P)\ge 0$)
determines a facet of the closed cone of \LV-vectors of Eulerian posets of 
rank $n+1$.
\end{proposition}
\begin{proof}
Apply the Facet Lemma~\ref{facetlemma} with $M=\emptyset$.
For a nonempty even set $R$, the interval system $\iv[R]$ of $R$ is nonempty,
so $\sum_{Q\subseteq[1,n]} \lv_Q(P(n,\iv[R]))
=\sum_{{\cal J}\subseteq \iv[R]} (-1)^{|{\cal J}|} =0$.
\end{proof}

\begin{theorem}
\label{facets}
Let $V$ be a subset of $[1,n]$ such that every
$I\in \iv[V]$ has 
cardinality at least $2$, and every $I\in \iv[[0,n+1]\setminus V]$ has
cardinality at most $3$. Assume that $M$ is a subset of $V$ such that
every $[a,b]\in\iv[V]$ satisfies the following: 
\begin{enumerate}
\item[{\em ($i$)}] $M\cap [a,b]= \emptyset,$ $[a,a+1]$, or $[b-1,b]$. 
\item[{\em ($ii$)}] If $a\not \in M$ then $a-2\in \{-1\}\cup M$.
\item[{\em ($iii$)}] If $b\not \in M$ then $b+2\in \{n+2\}\cup M$.
\end{enumerate}
Then 
\begin{eqnarray}
\label{facetineq}
(-1)^{|M|/2}\sum_{M\subseteq Q\subseteq V}L_Q(P)\ge 0
\end{eqnarray}
determines a facet of ${\cal C}^{n+1}_{\cal E}$.
Furthermore, if we strengthen {\em ($i$)} by also requiring 
$M\cap [a,a+2]=\emptyset$
for every $[a,a+2]\in\iv[V]$, then distinct pairs $(M,V)$ give distinct facets.
\end{theorem}

\begin{proof}
If $M=\emptyset$, then conditions ($ii$) and ($iii$) force $V=[1,n]$ (or 
$V=\emptyset$ if $n\le 1$).
The resulting inequality, $\sum_{Q\subseteq[1,n]} L_Q(P)\ge 0$, gives
a facet, as shown in Proposition~\ref{ineq0}.
Now assume that $M\neq\emptyset$.

Step 1 is to prove that inequality (\ref{facetineq}) holds for all Eulerian
posets.
Note that $\iv[M]$ is a nonempty collection of intervals of length two.  From
each such interval choose one endpoint adjacent to an element of
\mbox{$[0,n+1]\setminus V$}.
Let $T$ be the set of these chosen elements.
The Inequality Lemma \ref{ineqlemma} applies to these $T$ and $V$ 
because each interval of $V$ contains at most one interval of $\iv[M]$,
and hence at most one element of $T$.
The resulting inequality is 
$(-1)^{|T|}\sum_{T\subseteq Q\subseteq V} L_Q \ge 0$.
Now $L_Q(P)=0$ for all $P$ if $\iv[Q]$ contains an odd interval.
So we can restrict the sum to even sets $Q$.
Since $Q$ must be contained in $V$, such a $Q$ must contain the 
intervals of $M$.
Thus, 
$(-1)^{|M|/2}\sum_{M\subseteq Q\subseteq V}L_Q(P)\ge 0$.

Step 2 is to prove that if $I\subseteq [1,n]$ is an interval of
cardinality at least 2 and $I$ contains an element $i$ not in $V$,
then $I$ contains an element adjacent to an interval of $M$.
If an interval from $\iv[V]$ ends at $i-1$, then either $i-1\in M$ or 
$i+1\in M$ by ($iii$) (since $i+1<n+2$).
Similarly, if an interval from $\iv[V]$ begins at $i+1$, then either $i-1\in M$ or $i+1\in M$.
So assume no interval from $\iv[V]$ begins at $i-1$ or ends at $i+1$.
The hypothesis of the theorem states that every interval from 
$\iv[[0,n+1]\setminus V]$ has cardinality at most three. 
Thus the
interval $[i-1,i+1]$ belongs to
$\iv[[0,n+1]\setminus V]$. Hence $i-2\in \{-1\}\cup V$ and 
$i+2\in \{n+2\}\cup V$. If $i-2=-1$ then $I\supseteq
[i,i+1]=[1,2]$, condition ($ii$) applied to $a=3$ yields $3\in M$, and 
$2\in I$ is adjacent to $3$. The case when $i+2=n+2$ is dealt with
similarly. Finally, if $i-2$ and $i+2$ are both endpoints of 
intervals from $\iv[V]$, then, since $i\not\in M\cup\{-1,n+2\}$,
condition ($ii$) applied to $a=i+2$ and
condition ($iii$) applied to $b=i-2$ yield $i+2\in M$ and $i-2\in M$. 
Either $i-1$ or $i+1$ belongs to $I$ and each of them is adjacent to
an element of $M$. 

Recall that for $\iv$ an even interval system, the vector                       
$(\lv_Q(P(n,\iv)): \mbox{$Q\subseteq [1,n]$})$  is in the closed cone of 
$\lv$-vectors of half-Eulerian posets. 
Step 3 is to show that for each even set $R\ne M$, there exists an even         
interval system $\iv$ with $\cup_{i\in\iv}I=R$ such that 
$(-1)^{|M|/2}\sum_{M\subseteq Q\subseteq V}
\lv_Q(P(n,\iv))=0$.                                                       

Let $R$ be an even set not equal to $M$. 
If $M\not\subseteq R$, then for every $Q$ containing $M$,
$\lv_Q(P(n,\iv[R]))=0$. 
Now suppose $M\subseteq R$, but $R\not\subseteq V$.
Let $I$ be an interval of $\iv[R]$ such that $I\not\subseteq V$.
Then $I$ contains an element adjacent to an interval of $M$.
Since $M\subseteq R$ and $I$ is a maximal interval in $R$,                      
$I\cap M\ne\emptyset$.                                                          
Thus every union of intervals of $\iv[R]$ containing $M$ must contain $I$ and   
thus an element not in $V$.
So $\sum_{M\subseteq Q\subseteq V}\lv_Q(P(n,\iv[R]))=0$,       
because all terms are zero.

Finally, suppose $M\subseteq R\subseteq V$ and $R\ne M$.         
Let $\iv$ be the interval system of $R$ consisting only of intervals of         
length 2.                                                                       
Then every interval of $M$ is in $\iv$.                                         
This is because every interval of $M$ is of length 2, with at least one of
its endpoints adjacent to an element not in $V$.
So
$\sum_{M\subseteq Q\subseteq V} \lv_Q(P(n,\iv))
=\sum_{\iv[M]\subseteq {\cal J}\subseteq \iv} (-1)^{|{\cal J}|}=0$,
since $R\ne M$ implies $\iv\ne\iv[M]$.

By the Facet Lemma \ref{facetlemma}, the inequality
$(-1)^{|M|/2}\sum_{M\subseteq Q\subseteq V}L_Q(P)\ge 0$
gives a facet of ${\cal C}^{n+1}_{\cal E}$.

Now we show that under the added 
condition $M\cap [a,a+2]=\emptyset$ for every $[a,a+2]\in \iv[V]$, 
the facets obtained are distinct.

Note that two $(M,V)$ pairs can give the same inequality only if they have the
same $M$, because $\LV_M$ is included in the linear form for $(M,V)$, 
and $M$ is the minimal (by set inclusion) set for which $L_M$ is in the 
form.
Now for fixed $M$, we show that $(M,V_1)$ and $(M,V_2)$ give distinct linear 
inequalities when $V_1\neq V_2$. Since the sets $V_1$ and $V_2$ are
different, there is an interval $[a,b]$ such that 
$[a, b]$ occurs in exactly one of $\iv[V_1]$ or $\iv[V_2]$. Let $[a,b]$ 
be a maximal interval with this property. Without loss of generality 
assume $[a,b]\in\iv[V_1]$. 
Then $[a,b]$ is contained in no interval of $\iv[V_2]$.

{\em Case 1.} $M\cap [a,b]=\emptyset$. 
Then for every $i$, $a\le i\le b-1$, the term $L_{[i,i+1]\cup M}$ 
occurs in the inequality for $(M,V_1)$.
At least one of these terms does  not occur in the inequality for
$(M,V_2)$, because $[a,b]\not\subseteq V_2$.

{\em Case 2.} $M\cap [a,b]=[a,a+1]$.
Since $M\subseteq V_2$ and $[a,b]\not\subseteq V_2$, $b>a+1$.
By the strengthened hypothesis on $M$,  $b\geq a+3$. 
Then for every $i$, $a+2\le i\le b-1$, the term $L_{[i,i+1]\cup M}$ 
occurs in the inequality for $(M,V_1)$.
At least one of these terms does  not occur in the inequality for
$(M,V_2)$, because $[a,b]\not\subseteq V_2$.

{\em Case 3.} $M\cap [a,b]=[b-1,b]$.  The proof is similar to Case 2.

Thus, with the condition
$M\cap [a,a+2]=\emptyset$ for every $[a,a+2]\in \iv[V]$, the facets given by
the theorem are all distinct.
\end{proof}

Theorem \ref{facets} may be restated and interpreted in terms of the
convolution of chain operators. We refer the interested reader to
\ref{S_ring} for that approach.

With the aid of PORTA (\cite{porta}), we verified that the theorems above
give all the extremes and facets of the Eulerian cone for rank at most 6.

\begin{theorem}
\label{rank6}
For rank $n+1\le 6$, the closed cone ${\cal C}^{n+1}_{{\cal E}}$
of flag vectors of Eulerian posets is finitely generated.
It has ${n \choose \lfloor n/2\rfloor}$ 
extreme rays, all generated by the flag vectors of the limit posets 
$\HD P(n,\iv)$ for $\iv$ even interval systems on $[1,n]$.
It has ${n \choose \lfloor n/2\rfloor}$ facets, all 
given by Proposition~\ref{ineq0} and Theorem~\ref{facets}.
\end{theorem}

At rank 7 the situation changes for both extreme rays and facets.

\begin{theorem}
\label{rank7}
{\em ($i$)} The cone ${\cal C}^7_{\cal E}$ is finitely generated,
with 24 extreme rays.
Twenty of the extreme rays are generated by 
the flag vectors of the limit posets 
$\HD P(n,\iv)$ for $\iv$ even interval systems on $[1,6]$.

{\em ($ii$)} The cone ${\cal C}^7_{\cal E}$ has 23 facets. 
Fifteen of the facets are given by the inequalities of 
Theorem~\ref{facets}.
Four additional facets come from the Inequality
Lemma~\ref{ineqlemma}.
The remaining four come from Theorem~\ref{ijk-ineqs}.
\end{theorem}

The four special extreme rays of the rank 7 Eulerian cone have corresponding
rays in the half-Eulerian cone.
The generators for the half-Eulerian cone are all obtained by adding the
flag vectors of limit posets associated with noneven interval systems.
The  summands do not satisfy the conditions of Proposition~\ref{D_hE} for
half-Eulerian posets, but the sum does.
The calculations are easily done in terms of the $\lv$-vector, using
Proposition~\ref{P_limitl}.
Specific sequences of half-Eulerian posets have been constructed
whose flag vectors converge to these four extremes.
The half-Eulerian posets are obtained by ``gluing together'' posets 
for each summand.
These are then converted to Eulerian posets by the horizontal doubling 
operation.
Below are the sums of limit posets used.
Descriptions of the half-Eulerian posets are found in \ref{limit_posets}.

Extreme 1: $P\left(6,\{[1,2],[2,6]\}+\{[2,5],[5,6]\}\right)$

Extreme 2: $P\left(6,\{[1,3],[3,4],[4,6]\}+\{[1,2],[2,3]\}
+\{[4,5],[5,6]\}\right)$

Extreme 3: $P\left(6,\{[1,2],[3,4],[4,5]\}+\{[3,5],[5,6]\}
+\{[1,2],[2,5]\}\right)$

Extreme 4: $P\left(6,\{[1,2],[2,4]\}+\{[2,5],[5,6]\}
+\{[2,3],[3,4],[5,6]\}\right)$

\vspace{6pt}

Note that for rank at most 7, the two cones ${\cal C}^{n+1}_{\cal D}$ and 
${\cal C}^{n+1}_{\cal E}$ are equal, because the generators of extreme rays 
specified in Theorems~\ref{rank6} and~\ref{rank7} are horizontal doubles
of half-Eulerian limit posets.

Perhaps all the extreme rays of the half-Eulerian cone (if not the Eulerian
cone)
can be obtained by gluing together Billera-Hetyei limit posets.

A complete description of the closed cone of flag vectors of Eulerian
posets remains open, and, as mentioned before, the cone is not even known to be
finitely generated.
We do not know if convolutions of the inequalities of 
Proposition~\ref{ineqlemma} and 
Theorem~\ref{ijk-ineqs} completely determine the cone.
A better understanding of the construction of extreme rays as sums of 
Billera-Hetyei limit posets would be valuable.

The study of Eulerian posets is motivated in part by questions about convex 
polytopes.
Is the cone of flag vectors of all Eulerian posets the same as or close to the
cone of flag vectors of polytopes?
The answer is no.
The inequalities of Proposition~\ref{ineqlemma} can be strengthened considerably
for polytopes.
The proof of Proposition~\ref{ineqlemma} uses only the fact that in an Eulerian
poset each interval has at least two elements of each rank.
For convex polytopes, each interval is at least the size of a Boolean algebra
of the same rank.
Thus, for example, where Proposition~\ref{ineqlemma} gives that 
$f_{1479}(P)-2f_{179}(P)\ge 0$ for Eulerian posets, for convex polytopes
the inequality $f_{1479}(P)-20f_{179}(P)\ge 0$ holds, because the rank 6
Boolean algebra has ${6\choose 3}=20$ elements of rank 3.
For ranks 4 through 7, we have verified that none of the extreme rays of the
Eulerian cone is in the closed cone of flag vectors of convex polytopes.


\appendix

\renewcommand{\thesection}{Appendix \Alph{section}}

\section{Some half-Eulerian limit posets of rank~$7$}

\label{limit_posets}

\renewcommand{\thesection}{\Alph{section}}
\renewcommand{\thesubsection}{\thesection.\arabic{subsection}}

Here are the constructions of half-Eulerian posets whose doubles give
Extremes~1, 2 and~3 of ${\cal C}^7_{\cal E}$.
Extreme~4 is the dual of Extreme~3.

In the following, $C^7$ denotes a chain of rank $7$.

\subsection{$P\left(6,\{[1,2],[2,6] \}+\{[2,5],[5,6] \}\right)$}

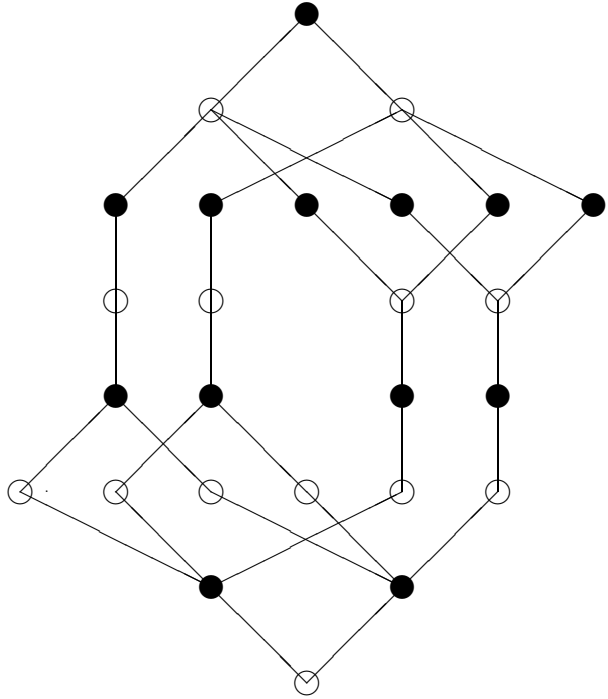
\begin{figure}[h]
\setlength{\unitlength}{3947sp}
\begingroup\makeatletter\ifx\SetFigFont\undefined
\gdef\SetFigFont#1#2#3#4#5{
  \reset@font\fontsize{#1}{#2pt}
  \fontfamily{#3}\fontseries{#4}\fontshape{#5}
  \selectfont}
\fi\endgroup
\begin{center}
\begin{picture}(3766,4366)(1118,-4644)
\thinlines
\put(3001,-4561){\circle{150}}
\put(1801,-3361){\circle{150}}
\put(1201,-3361){\circle{150}}
\put(2401,-3361){\circle{150}}
\put(3001,-3361){\circle{150}}
\put(1801,-2161){\circle{150}}
\put(2401,-2161){\circle{150}}
\put(2401,-961){\circle{150}}
\put(3601,-961){\circle{150}}
\put(3601,-3961){\circle*{150}}
\put(2401,-3961){\circle*{150}}
\put(1801,-2761){\circle*{150}}
\put(2401,-2761){\circle*{150}}
\put(1801,-1561){\circle*{150}}
\put(2401,-1561){\circle*{150}}
\put(3001,-361){\circle*{150}}
\put(3001,-1561){\circle*{150}}
\put(3601,-1561){\circle*{150}}
\put(4201,-1561){\circle*{150}}
\put(4801,-1561){\circle*{150}}
\put(3601,-3361){\circle{150}}
\put(4201,-3361){\circle{150}}
\put(3601,-2161){\circle{150}}
\put(4201,-2161){\circle{150}}
\put(3601,-2761){\circle*{150}}
\put(4201,-2761){\circle*{150}}
\put(1276,-3361){\makebox(1.6667,11.6667){\SetFigFont{5}{6}{\rmdefault}{\mddefault}{\updefault}.}}
\put(2401,-3961){\line( 1,-1){600}}
\put(3001,-4561){\line( 1, 1){600}}
\put(1201,-3361){\line( 2,-1){1200}}
\put(2401,-3961){\line(-1, 1){600}}
\put(2401,-3361){\line( 2,-1){1200}}
\put(3601,-3961){\line(-1, 1){600}}
\put(1201,-3361){\line( 1, 1){600}}
\put(1801,-3361){\line( 1, 1){600}}
\put(1801,-2761){\line( 1,-1){600}}
\put(2401,-2761){\line( 1,-1){600}}
\put(1801,-2761){\line( 0, 1){1200}}
\put(2401,-2761){\line( 0, 1){1200}}
\put(1801,-1561){\line( 1, 1){600}}
\put(2401,-1561){\line( 2, 1){1200}}
\put(2401,-961){\line( 1, 1){600}}
\put(3001,-361){\line( 1,-1){600}}
\put(2401,-961){\line( 1,-1){600}}
\put(2401,-961){\line( 2,-1){1200}}
\put(3601,-961){\line( 1,-1){600}}
\put(3601,-961){\line( 2,-1){1200}}
\put(3001,-1561){\line( 1,-1){600}}
\put(3601,-1561){\line( 1,-1){600}}
\put(4201,-1561){\line(-1,-1){600}}
\put(4801,-1561){\line(-1,-1){600}}
\put(3601,-2161){\line( 0,-1){1200}}
\put(4201,-2161){\line( 0,-1){1200}}
\put(2401,-3961){\line( 2, 1){1200}}
\put(3601,-3961){\line( 1, 1){600}}
\end{picture}
\end{center}
\caption{$P\left(6,\{[1,2],[2,6] \}+\{[2,5],[5,6] \}\right)$}
\label{F_lpos1}
\end{figure}

Take $\HD^{N}_{[1,2]}\HD^{N}_{[2,6]}\left(C^7\right)$ and 
$\HD^{N}_{[1,5]}\HD^{N}_{[5,6]}\left(C^7\right)$. Identify the elements of both
posets at rank $1$ and at rank $6$. Figure \ref{F_lpos1} represents the
resulting poset for $N=2$.

\newpage

\subsection{$P\left(6,\{[1,3],[3,4],[4,6]\}
+\{[1,2],[2,3]\}
+\{[4,5],[5,6] \}\right)$}

\begin{figure}[h]
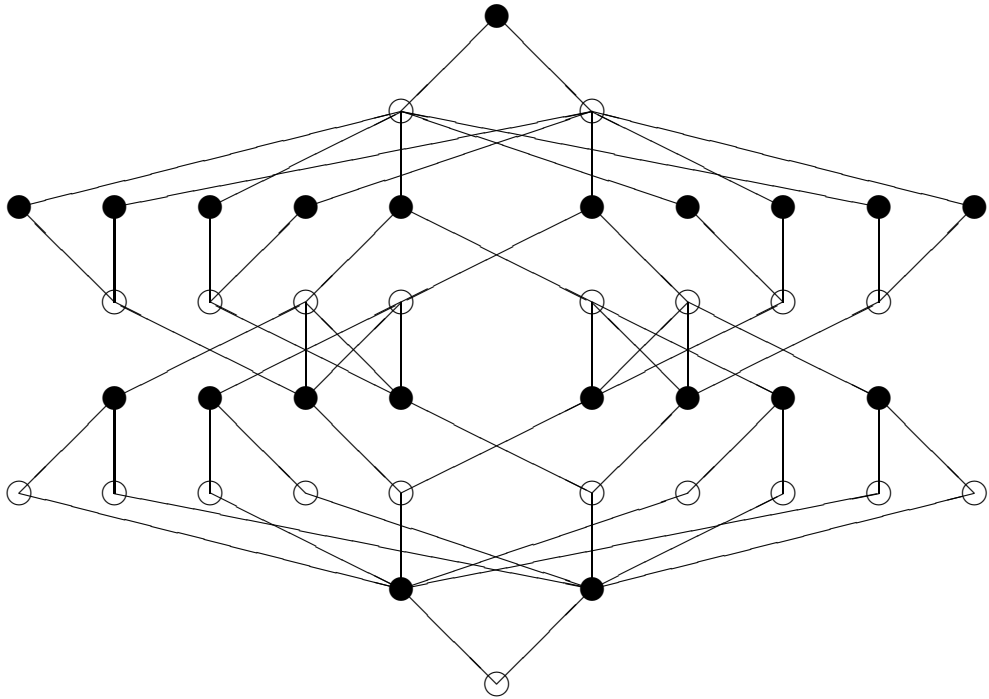

\setlength{\unitlength}{3947sp}
\begingroup\makeatletter\ifx\SetFigFont\undefined
\gdef\SetFigFont#1#2#3#4#5{
  \reset@font\fontsize{#1}{#2pt}
  \fontfamily{#3}\fontseries{#4}\fontshape{#5}
  \selectfont}
\fi\endgroup
\begin{center}

\end{center}
\caption{$P\left(6,\{[1,3],[3,4],[4,6]\}
+\{[1,2],[2,3]\}
+\{[4,5],[5,6] \}\right)$}
\label{F_lpos2}
\end{figure}

Take
\begin{eqnarray*}
P^I(N)&=&
\HD^{N}_{[1,3]}\HD^{N}_{[3,4]}\HD^{N}_{[4,6]}\HD^{N+1}_{[4,5]}
(C^7)\\ \\
P^{II}(N)&=&
\HD^{N+1}_{[1,2]}\HD^{N}_{[1,6]}\HD^{N}_{[2,4]}(C^7),
\quad\mbox{and}\\   \\
P^{III}(N)&=&
\HD^{N}_{[1,5]}\HD^{N}_{[3,5]}\HD^{N}_{[5,6]}(C^7).\\
\end{eqnarray*}
Identify the elements of $P^{I}(N)$ with the elements of $P^{II}(N)$ at
ranks $1,4,5$, and $6$. 
Identify the elements of $P^{I}(N)$ with the elements of $P^{III}(N)$ at
ranks $1,2,3$, and $6$. Figure \ref{F_lpos2} represents the resulting poset for
$N=2$.

\newpage

\subsection{$P\left(6,\{[1,2],[3,4],[4,5]\}
+\{[3,5],[5,6]\}
+\{[1,2],[2,5] \}\right)$}

Take

\begin{figure}
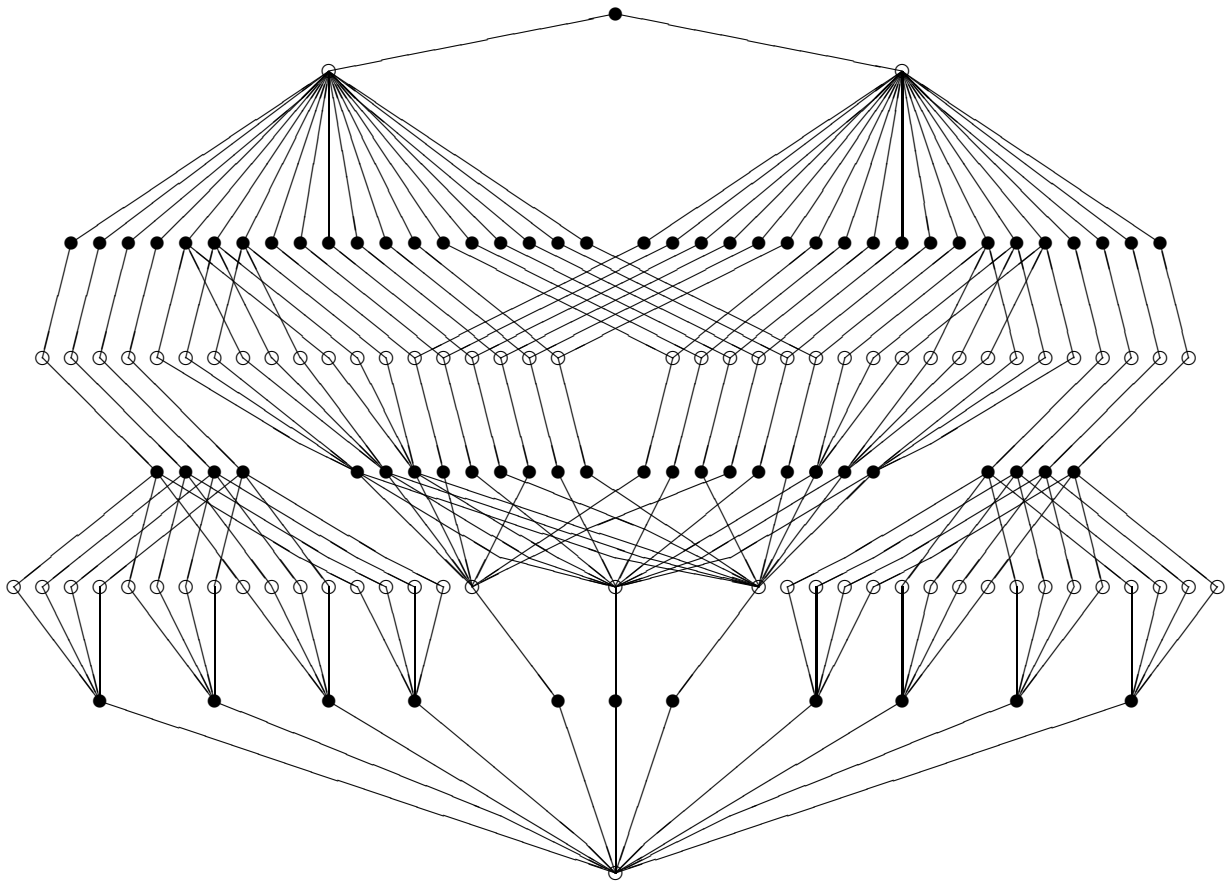

\setlength{\unitlength}{2368sp}
\begingroup\makeatletter\ifx\SetFigFont\undefined
\gdef\SetFigFont#1#2#3#4#5{
  \reset@font\fontsize{#1}{#2pt}
  \fontfamily{#3}\fontseries{#4}\fontshape{#5}
  \selectfont}
\fi\endgroup
\hspace*{-.75in}
%
\caption{$P\left(6,\{[1,2],[3,4],[4,5]\}
+\{[3,5],[5,6]\}
+\{[1,2],[2,5] \}\right)$}
\label{F_lpos3}
\end{figure}

\vspace*{-24pt}

$$
\begin{array}{rcll}

\vspace*{6pt}

P^I(N)&=&
\HD^{N+1}_{[1,2]}\HD^{N+1}_{[3,4]}\HD^{N}_{[3,6]}\HD^{N+1}_{[4,5]}
(C^7)&\mbox{(Figure \ref{F_lpos3a})}\\  

\vspace*{6pt}

P^{II}(N)&=&
\HD^{N+1}_{[1,5]}\HD^{N^2}_{[3,5]}\HD^{N}_{[5,6]}(C^7)
&\mbox{(Figure \ref{F_lpos3b}), and}\\  
P^{III}(N)&=&
\HD^{N+2}_{[1,2]}\HD^{N^2-N+2}_{[2,5]}\HD^{N}_{[1,6]}(C^7)&
\mbox{(Figure \ref{F_lpos3c}).}
\end{array}
$$
Identify the elements of $P^{I}(N)$ with the elements of $P^{II}(N)$ at
ranks $1,2$, and $6$. 
Identify the elements of $P^{I}(N)$ with the elements of $P^{III}(N)$ at
rank $6$. Figure \ref{F_lpos3} represents the resulting poset for
$N=2$.

\begin{figure}[h]
\setlength{\unitlength}{2368sp}
\begingroup\makeatletter\ifx\SetFigFont\undefined
\gdef\SetFigFont#1#2#3#4#5{
  \reset@font\fontsize{#1}{#2pt}
  \fontfamily{#3}\fontseries{#4}\fontshape{#5}
  \selectfont}
\fi\endgroup
\begin{center}
%
\end{center}
\caption{$P^I(2)$}
\label{F_lpos3a}
\end{figure}

\begin{figure}
\setlength{\unitlength}{2368sp}
\begingroup\makeatletter\ifx\SetFigFont\undefined
\gdef\SetFigFont#1#2#3#4#5{
  \reset@font\fontsize{#1}{#2pt}
  \fontfamily{#3}\fontseries{#4}\fontshape{#5}
  \selectfont}
\fi\endgroup
\begin{center}
%
\end{center}
\caption{$P^{II}(2)$}
\label{F_lpos3b}
\end{figure}

\begin{figure}
\setlength{\unitlength}{2368sp}
\begingroup\makeatletter\ifx\SetFigFont\undefined
\gdef\SetFigFont#1#2#3#4#5{
  \reset@font\fontsize{#1}{#2pt}
  \fontfamily{#3}\fontseries{#4}\fontshape{#5}
  \selectfont}
\fi\endgroup
\hspace*{-.75in}
%
\hspace*{-.25in}
\caption{$P^{III}(2)$}
\label{F_lpos3c}
\end{figure}

\pagebreak

\renewcommand{\thesection}{Appendix \Alph{section}}
\section{The Billera-Liu ring of chain operators}
\label{S_ring}

\renewcommand{\thesection}{\Alph{section}}

As in Billera and Liu (\cite{Billera-Liu}) we view the flag $f$-vector as
a vector of {\em chain operators} $\left(f_S^{n+1}\::\: S\subseteq
[1,n]\right)$; here $f_S^{n+1} (P)= f_S(P)$ if $P$ is a graded poset of
rank $n+1$ and $0$ otherwise.
The following multiplication of chain operators $f^n_S$ ($n\geq 1$,
$S\subseteq [1,n-1]$) was introduced by Kalai
in \cite{Kalai} and studied for Eulerian posets by Billera and Liu in
\cite{Billera-Liu}: 
$$f^m_S f^n_T\df f^{m+n}_{S\cup \{m\}\cup (T+m)}.$$
It is straightforward that given a pair of valid linear inequalities 
$$F=\sum_{S\subseteq [1,m-1]} a_S f^{m}_S \geq 0 \quad\mbox{and}\quad
G=\sum_{T\subseteq [1,n-1]} b_S f^{n}_S \geq 0$$ 
that hold for a class of graded posets, the linear inequality 
$FG\geq 0$ is also valid for the same class. It was
observed by Billera and Liu in \cite[Proposition 1.3]{Billera-Liu} that
for the class of all graded posets the converse holds as well: if
$FG\geq 0$ is a valid inequality, then either both $F\geq 0$ and
$G\geq 0$ are valid inequalities, or both  $-F\geq 0$ and $-G\geq 0$ are valid
inequalities. According to \cite[Theorem 2.1]{Billera-Liu} the
associative algebra generated by all chain operators (whose domain is
taken to be the class of all graded posets) is the free polynomial ring
in variables $\{f^i_{\emptyset}\::\: i\geq 1\}$. If we take the degree 
of the variable $f^i_{\emptyset}$ to be $i$, then linear combinations of the
form $F=\sum_{S\subseteq [1,m-1]} a_S f^{m}_S$ become homogeneous
polynomials. Hence, as noted by Billera and Hetyei in
\cite{Billera-Hetyei}, one can use a result of Cohn in \cite[Theorem
3]{Cohn} that the semigroup of homogeneous polynomials of a free
graded associative algebra has unique factorization. Hence an inequality
can be checked factor-by-factor.
Billera and Hetyei also showed in \cite{Billera-Hetyei} that
for the class of all graded posets the product of two facet inequalities
is almost always a facet inequality, every exception being a consequence 
of the equalities  
$$
f^m_{\emptyset}f^n_{\emptyset}=f^{m+n}_m=
\left(f^{m+n}_m-f^{m+n}_{\emptyset}\right)+f^{m+n}_{\emptyset}.
$$  

For Eulerian and half-Eulerian posets, it is advisable to convert our
expressions into the flag-$\lv$ or flag-$L$ forms respectively. 
Straighforward substitution into the definition shows 
$$\lv_S^m\lv_T^n=\lv_{S\cup (T+m)}^{m+n} \quad\mbox{and}
\quad \LV_S^m\LV_T^n=2 \LV_{S\cup (T+m)}^{m+n}$$
This means that when we write $[u_S]=\LV^{n}_S$ as the coefficient of 
the $ce$-word $u_S$, the convolution of the forms 
$\sum_{S\subseteq [1,m-1]} a_S [u_S]$ and
$\sum_{T\subseteq [1,n-1]} b_T [u_T]$
is a constant multiple of the form
$\sum_{S\subseteq [1,m-1]}\sum_{T\subseteq [1,n-1]} a_S b_T [u_S c
u_T].$
In particular, if only monomials of $c$ and $ee$ occur in each factor,
the same holds for the convolution. Hence the same
result of Cohn \cite[Theorem 3]{Cohn} on unique homogeneous
factorization proves the following. 

\begin{proposition}
Every homogeneous linear form $\sum_{S\subseteq [1,n]} a_S \lv^{n+1}_S$ 
or\\
$\sum_{S\subseteq [1,n]} a_S \LV^{n+1}_S$, where $S$ ranges over 
only even sets, can be uniquely written as a product of irreducible
expressions of the same kind.
\end{proposition}
Let us call such expressions {\em even $\lv$-forms} and {\em even
$\LV$-forms}, respectively. The interest in this factorization stems from 
the following observation.

\begin{proposition} 
Let $F$ and $G$ both be even $\lv$-forms. Then $FG\geq 0$ holds for all
half-Eulerian posets if and only if either both $F\geq 0$ and $G\geq 0$
or both  $-F\geq 0$ and $-G\geq 0$ hold for all half-Eulerian
posets. The analogous statement is true for even $\LV$-forms and
Eulerian posets. 
\end{proposition}
Only the ``only if'' implication is not completely trivial. In the
half-Eulerian case, all we need to observe is that for a pair $(P,Q)$
of half-Eulerian posets the poset $P\circ Q$ obtained by putting all elements
of $Q$ above all elements of $P$, and identifying the top element of $P$
with the bottom element of $Q$, is half-Eulerian. Moreover, if for
posets $P_1,P_2$, and $Q$ and forms $F$ and $G$, $F(P_1)>0$,
$F(P_2)<0$, and $G(Q)>0$, then $FG(P_1\circ Q)=F(P_1)
G(Q)>0$ and $FG(P_2\circ Q)=F(P_2)
G(Q)<0$. The same argument works for Eulerian posets using
$\HD^2_{\{\rank(P)\}} (P\circ Q)$ instead of $P\circ Q$.

In terms of convolutions, Proposition~\ref{ineqlemma} states that the product
of valid inequalities of the form $f^n_{\emptyset}\geq 0$ and
$f^n_{i}-2f^n_{\emptyset}\geq 0$ is a valid inequality for all Eulerian
posets. Theorem~\ref{facets} describes a subclass of these products that
yield facet inequalities. Using ideas extracted from the proof, one can
show the following, somewhat strengthened statements.

\begin{proposition}
If $F\geq 0$ defines a facet of ${\cal C}^{n+1}_{\cal E}$, then 
$F(f^{k+1}_1-2 f^{k+1}_{\emptyset})\geq 0$ defines a facet of
${\cal C}^{n+k+2}_{\cal E}$. 
\end{proposition}

\begin{proposition}
If $F\geq 0$ defines a facet of ${\cal C}^{n+1}_{\cal E}$, and 
$F$ can be written as 
$$F=\sum_{S\subseteq [1,n]} a_S  L^{n+1}_{S}$$
where $S$ ranges over only even sets that contain $n$, then 
$Ff^{k+1}_{\emptyset}\geq 0$ and 
$Ff^1_{\emptyset}f^1_{\emptyset}\geq 0$  
define facets of
${\cal C}^{n+k+2}_{\cal E}$ and ${\cal C}^{n+3}_{\cal E}$, respectively. 
\end{proposition}

It seems to be difficult, however, even in the case of these simple
factors to predict which products yield facet inequalities.
For example $(f^5_1-2f^5_{\emptyset})
f^1_{\emptyset}=(f^6_1-2f^6_{\emptyset})+{1\over 2}
(f^3_1-2f^3_{\emptyset})(f^3_1-2f^3_{\emptyset})\geq 0$ does not
define a facet of ${\cal C}^{6}_{\cal E}$, while it can be shown 
that $(f^5_1-2f^5_{\emptyset})f^3_{\emptyset}\geq 0$ defines a facet of
${\cal C}^{8}_{\cal E}$.

\end{document}